\documentclass[fleqn]{article}
\usepackage{amsmath}
\usepackage[]{url}
\usepackage[dvips]{graphicx}
\usepackage{authblk}
\usepackage{amssymb, amsmath, amsfonts, latexsym}
\usepackage{amssymb}

\newcommand{\qed}{\hfill \rule{1.6mm}{1.6mm}}
\newcommand{\no}{\nonumber}
\newcommand{\noi}{\noindent}

\catcode`\@=11
\def\theequation{\@arabic{\c@section}.\@arabic{\c@equation}}
\catcode`\@=12
\numberwithin{equation}{section}
\title{Game-theoretic approach to risk-sensitive benchmarked asset management. }

\author[*]{ Amogh Deshpande}
\author[**]{Saul Jacka}
\affil[*]{Institute for Financial and Actuarial Mathematics, Department of Mathematical Sciences, University of Liverpool and Department of Statistics, University of Warwick, UK. Email: addeshpa@gmail.com}
\affil[**]{Department of Statistics, University of Warwick ,UK. Email: S.D.Jacka@warwick.ac.uk}
\date{}

\begin{document}
\maketitle \baselineskip20pt
\parskip10pt
\parindent.4in
\begin{abstract} \noi  In this article we consider a game theoretic approach to the Risk-Sensitive Benchmarked Asset Management  problem  (RSBAM) of Davis and Lleo \cite{DL}. In particular, we consider a stochastic differential game between two players, namely, the investor who has a power utility while the second player represents the market which tries to minimize the expected payoff of the investor. The market does this by modulating a stochastic benchmark that the investor needs to outperform. We obtain an explicit expression for the optimal pair of strategies as for both the  players.
\end{abstract}
\noi {\bf Key Words}: Risk- Sensitive control ,  zero sum stochastic differential game.

\newpage

\section{Introduction }
$\indent$ In this article we shall develop a game theoretic version of a continuous time optimization model with risk-sensitive control approach more specifically termed as Risk-sensitive control portfolio optimization  (RSCPO). The RSCPO balances an investor's interest in maximizing the expected growth rate of wealth against his
aversion to risk due to deviations of the  realized rate from the expectation. The subjective notion of investor's risk aversion is parameterized
by a single variable, say $\theta$. More formally, we write the finite horizon risk-sensitive optimization criterion as : Maximize,\\
\begin{eqnarray*}
J_{T,h}:=-\frac{1}{\theta}\log{E[e^{-\theta F({T,h})}]}
\end{eqnarray*}
where $F(T,h)$ is the time-$T$ value reward function corresponding to control $h$. In the optimal investment problem we take $F(T,h)=\log{V(T)}$ where $V(t)$ is the time $t$-value of the portfolio  corresponding to portfolio asset allocation $h$.
{
An asymptotic expansion around $\theta=0$ for the above criterion yields
\begin{eqnarray*}
J_{T,h}=E[F({T,h})]-\frac{\theta}{2}Var(F({T,h}))+O(\theta^{2})
\end{eqnarray*}
From this expression it is clear this criterion compromises between maximizing the  portfolio return while penalizing the riskiness . The optimal expected utility function depends on $\theta$ and is a generalization of the traditional  stochastic control approach to utility optimization in the sense that now the degree of risk aversion  of the investor is explicitly parameterized through $\theta$ rather than importing it in the problem via an exogenous utility function. Values of $\theta>0$ correspond to a risk-averse investor, $\theta<0$ to a risk-seeking investor and $\theta=0$ to a risk-neutral investor  who maximizes}
\begin{eqnarray*}
J_{T,h}:={E[{ F(T,h)}]}
\end{eqnarray*}
There has been a substantial amount of research on the infinite-time horizon ergodic problem:
 \begin{eqnarray*}
&&  \max~~~ \bar{J}_{\infty}~~ \mbox{where}\\~~~~~
&&\bar{J}_{\infty}=\liminf_{t \rightarrow \infty}-\frac{1}{\theta}t^{-1}\log{E[e^{-\theta F(t,h)}]}
 \end{eqnarray*}
 Though these type of problems are interesting in their own right, they  are not readily applicable to practical asset management because of non-uniqueness of optimal controls.
 \\
 \indent In the past decade,  applications of risk-sensitive control to asset management have proliferated. Risk-sensitive control was first applied to solve financial problems by Lefebvre and Montulet \cite{LM} in a corporate finance context.  Fleming \cite{Fl} was the first to show that some investment optimization models could be reformulated as  risk-sensitive control problems. Bielecki and Pliska \cite{BP} considered a model with $n$ securities and $m$ economic factors with no transaction cost. They were the first to apply continuous-time risk-sensitive control as a practical tool that could be used to solve ``real-world'' portfolio selection problems. They considered a long-term asset allocation problem and proposed the logarithm of the investor's wealth as a reward function, so that the investor's objective is to maximize the risk-sensitive (log) return of his/her portfolio. They derived the optimal control and solved the associated Hamilton-Jacobi-Bellman (HJB) PDE under the restrictive assumption that  the securities and economic factors have independent noise. In \cite{BP03}, Bielecki and Pliska  went on to study the economic properties of the risk-sensitive asset management criterion and then extended the asset management model into an intertemporal CAPM in \cite{BP04}.
 Fleming and Sheu \cite{Fs} analyzed an investment model  similar to that of Bielecki and Pliska \cite{BP}. In their model, however, the factor process and the security price process were assumed correlated.
  A major contribution was made by Kuroda and Nagai \cite{KN} who introduced an elegant solution method based on a change of measure argument which transforms the risk sensitive control problem into a linear exponential of a quadratic regulator. They solved the associated HJB PDE over a finite time horizon and then studied the properties of the ergodic HJB PDE related to $\bar{J}_{\infty}$. Recently, Davis and Lleo \cite{DL} applied this change of measure technique to solve, for both the finite and an infinite horizon, a  risk-sensitive benchmark investment problem (RSBAM) in which an investor selects an asset allocation to outperform a given financial benchmark. In the Kuroda and Nagai set-up $\theta$ represents the sensitivity of an investor to total risk, whereas in the RSBAM, $\theta$  represents the investor�s sensitivity to active risk i.e.   additional risk the investor is willing to take in order to outperform the benchmark.
 It is obvious that for outperforming a stochastic benchmark, an investor will have to modify his or her optimal trading strategy. Then the question of  interest to us is: ``What is the investor's worst case strategy for an opposing stochastic benchmark"?. In particular, one can even take the jaundiced point of view that the benchmark  will be set retrospective to the worst case. For example, if a portfolio fund manager outperforms the set benchmark, the principal may remark this out-performance either as best achieved or poorly achieved with respect to the underlying worst-case scenario. So, in this article we consider a game-theoretic version of the problem  within the benchmark framework of Davis and Lleo \cite{DL}. In it, we consider a stochastic differential game between two players, namely, the investor (who has a power utility) and a second player, representing the market, who tries to minimize the expected payoff of the investor. We explicitly characterize the optimal allocation of assets and the optimal choice of benchmark index. \\
 \indent In this article, we consider the benchmark process ex-ante that evolves according to a controlled diffusion process. We contrast this approach to the one  of Heath and Platen   \cite{HP}. In their methodology, they use the growth optimal portfolio itself as a benchmark which is closer to the concept of the numeraire portfolio. Although there has been a long history of applying risk-sensitive optimal control to problems in finance, a game-theoretic version of such problems in finite horizon is missing from the literature. We intend to elaborate further on this now.\\
 \indent In the next section we briefly describe the framework of the risk-sensitive zero sum stochastic differential game corresponding to the desired game ({\bf P1})( refer 2.8a). In the third section we reformulate the objective criterion under evaluation as a linear exponential of quadratic regulator problem ({\bf P2}) (refer 3.11). In the fourth section we provide a verification lemma that will help us solve this game problem. In the fifth section we derive the optimal controls and obtain an explicit expression for the associated value of the game. The article as usual concludes with remarks and pointers to future direction of work.\\
\indent Broadly speaking our aim is to derive the saddle-point equilibrium pair for the game ({\bf P1}). To achieve this, we first obtain saddle point strategy for the game ({\bf P2}). We then show that  the saddle point equilibrium for ({\bf P2}) is also saddle point equilibrium for ({\bf P1}). \\
\section{Risk-sensitive zero sum stochastic differential game}
We consider a market consisting of $m+1 \geq 2$ securities with $n\geq 1$ factors. We assume that the set of securities includes one bond whose price is governed  by the ODE
\begin{eqnarray}\label{2.1}
dS_{t}^{0}=r_{t}S_{t}^{0}dt,~~~~ S_{0}^{0}=s^{0}
\end{eqnarray}
where $r_{t}$ is a deterministic function of $t$. The other security prices and factors are assumed to satisfy the following SDE's
\begin{eqnarray}\label{2.2}
dS_{t}^{i}=S_{t}^{i}\{(a+AX_{t})^{i}dt+\sum_{k=1}^{n+m}{\sigma_{k}^{i}dW_{t}^{k}}\}, ~~~~S_{0}^{i}=s^{i}, i=1,...,m,
\end{eqnarray}
where the factor process $X_{t}$ satisfies,
\begin{eqnarray}\label{2.3}
dX_{t}=\{(b+BX_{t})dt+\Lambda dW_{t}\}, X_{0}=x \in \mathbb{R}^{n}
\end{eqnarray}
Here $W_{t}=(W_{t})_{k=1,...,n+m}$ is an $n+m$ dimensional standard Brownian motion  defined on a filtered probability space ($\Omega,\mathcal{F},\mathbb{P},\mathcal{F}_{t}$). \\ \indent The factor process can represent  macro-economic indicators such as GDP, inflation and market index data. The  stock price dynamics are modulated by the factor process. Hence one can incorporate the effect of macro-economic indicators into the investment optimization problem by using the stock price process modulated by the factor process $X_t$. \\
\indent The model parameters $A, B, \Lambda$ are respectively $m \times n, n \times n, n \times (m+n)$ constant matrices and $a \in \mathbb{R}^{m}$, $b \in \mathbb{R}^{n}$. The constant matrix $(\sigma_{k}^{i})_{\{i=1,2....,m; k=1,2,...,(n+m)}\}$ will be denoted by $\Sigma$ in what follows.   \\ \indent In Kuroda and Nagai \cite{KN} it is assumed that the factor process and the stock price process do not have independent noise i.e. $\Sigma\Lambda^{'}\neq 0$. This assumption is in sharp contrast to Bielecki and Pliska \cite {BP} who conversely assume that $\Sigma\Lambda^{'}=0$. We will assume that $\Sigma\Lambda^{'}\neq 0$.  \\
\indent Let $\mathcal{G}_{t}=\sigma(S_{u},X_{u},L^{\gamma}_{u};u \leq t)$ be the sigma-field generated by the underlying stock price process, factor process and benchmark process $L^{\gamma}$ to be defined later up to time $t$. The  investment strategy which represents the proportional allocation of total wealth in the $i^{th}$ security $S_{t}^{i}$ is denoted by $h_{t}^{i}$  for $i=1,...,m$. Strategy $(h^{0}_{t},h_{t})_{0 \leq t \leq T}$ is said to be an investment strategy up to time $T$. We set $S^{'}_{t}:=(S_{t}^{1},S_{t}^{2},...,S_{t}^{m})^{'},  h^{'}_{t}:=(h_{t}^{1},...,h_{t}^{m})^{'}$. The space of controls ${\mathcal{H}}(T)$ consists of $\mathbb{R}^{m}$-valued controls for the investor as follows: ${\mathcal{H}}(T)$   is the set of $\{\mathcal{B}[0,T] \otimes \mathcal{G}_{t}\}_{\{t \geq 0\}}$-progressively measurable stochastic processes such that $\sum_{i=1}^{m}{h_{t}^{i}}+h_{t}^{0}=1$ and where $P(\int_{0}^{T}{|h_{s}|^{2}}ds<\infty)=1~~\forall ~T<\infty $ and ${E[e^{\int_{0}^{T}{\theta^{2}h^{'}_{s}\Sigma\Sigma^{'}h_s ds}}]}^{\frac{1}{2}}< \infty$.\\
\indent For given $h \in {\mathcal{H}}(T)$, the process $V_{t}=V_{t}^{h}$ represents the investor's wealth at time $t$, under the control $h$, and satisfies the following SDE dynamics,
\begin{eqnarray*}
\frac{dV^{h}_{t}}{V^{h}_{t}}&=&(r_{t} + h_{t}^{'}((a+AX_{t})-r_{t}1))dt+ h_{t}^{'}\Sigma dW_{t}; V^{h}_{0}=v
\end{eqnarray*}
which can be rewritten  as,
\begin{eqnarray}\label{2.4}
\frac{dV^{h}_{t}}{V^{h}_{t}}&=&(r_{t} + h_{t}^{'}d_{t})dt+ h_{t}^{'}\Sigma dW_{t}; V^{h}_{0}=v
\end{eqnarray}
where $d_{t} \triangleq a+AX_{t}-r_{t}1$. From  equation (\ref{2.4}) it can be seen that if $a+AX_{t}=r_{t}1$ i.e. $d_{t}=0$, then the portfolio wealth process evolves with drift equal to the riskless interest rate $r_{t}$.  We make an assumption here that the securities price volatility matrix $\Sigma$ is a full rank matrix. If it is not full-rank then $h^{'}\Sigma=0$ for some $h \neq 0$. Hence the market contains redundant asset(s) and the portfolio value process $V^{h}_{t}$ will grow at a rate different than the risk-less interest rate $r_{t}$ when $h^{'}d \neq 0$ resulting in an arbitrage. This is the case if the portfolio contains two or more redundant assets for example a stock and an option on the same stock. Hence we remove redundancy till
the resultant matrix $\Sigma$ is of full rank thereby ensuring that there exist  no further possibility of arbitrage by trading in the resultant portfolio. In our benchmark model we  express the objective through a new optimization criterion corresponding to a reward function $F$ which represents the log excess return of the asset portfolio over its benchmark and is given as
\begin{eqnarray*}
F(t;h,\gamma)=\log{\frac{V^{h}_{t}}{L^{\gamma}_{t}}}~~~F(0;h,\gamma)=\log{f}
\end{eqnarray*}
We now formally state the Risk-sensitive Benchmarked Asset management problem (RSBAM) that we  solve. \goodbreak
{\bf Problem : Risk-sensitive Benchmarked Asset Management (RSBAM)}\\
We first define the objective criterion $J$ as,\\
\begin{eqnarray}\label{2.5}
J(f,x,h,\gamma;T)& \triangleq & \frac{-2}{\theta}\log E[\exp{[\frac{-\theta}{2}{F(T,h,\gamma)}]}]\no\\
&=&\frac{-2}{\theta}\log{E[\bigg(\frac{V^{h}_{T}}{L^{\gamma}_{T}}\bigg)^{-\theta/2}]} \no \\
&=&\frac{-2}{\theta}\log{E[U(\frac{V^{h}_{T}}{L^{\gamma}_{T}})]}
\end{eqnarray}
where the utility function $U(\cdot)$ is  $U:x \rightarrow x^{{-\frac{\theta}{2}}}$. The dynamics of the benchmark process is a diffusion process $L^{\gamma}$ modulated by  a (Markovian) control $\gamma$ given by
\begin{eqnarray}\label{2.6}
\frac{dL_{t}^{\gamma}}{L_{t}^{\gamma}}=(\alpha_{t}+\beta_{t}X_{t})dt+\gamma^{'}_{t}dW_{t}
\end{eqnarray}
 where $\alpha_t \in \mathbb{R}$ and $\beta \in \mathbb{R}^{1 \times n}$.
The space of controls $\Gamma(T)$ consists of the market control represented by $\gamma$ that is $\mathbb{R}^{n+m}$-valued.  $\Gamma(T) $ consists of  progressively measurable controls measurable w.r.t to ${\{\mathcal{B}[0,T] \otimes \mathcal{G}_{t}\}}_{t\geq 0}$ and where $P(\int_{0}^{T}{|\gamma_{s}|^{2}}ds<\infty)=1~~\forall ~T<\infty$ and ${E[e^{\theta^{2}\int_{0}^{T}{\gamma^{'}_{s}\gamma_{s}d{s}}}]}^{\frac{1}{2}}< \infty$ .\\
By a simple application of Ito's formula we have:
\begin{eqnarray}\label{2.7}
dF(t,h,\gamma)=d\log(\frac{V^{h}_{t}}{L_{t}^{\gamma}})&=&F(t,h,\gamma)\{[{r_{t}+h_{t}{'}(a+AX_{t}-r_{t}1)-(\alpha_{t}+\beta_{t}X_{t})-\frac{1}{2}h_{t}^{'}\Sigma\Sigma^{'}h_{t}+\frac{1}{2}\gamma_{t}^{'}\gamma_{t}}]dt\no\\
&+&(h^{'}_{t}\Sigma-\gamma^{'}_{t})dW_{t}\}
\end{eqnarray}
\indent We are now in a position to formally state the game-theoretic version of the game. For a given $\theta>0$, we consider a stochastic differential game between  two players, namely, the investor (who has a power utility) $U$ and who modulates the payoff for given $\gamma \in \Gamma(T)$ via control  $h \in {\mathcal{H}}(T)$. On the other hand the second player, say the market,  behaves antagonistically to the investor by setting a benchmark for the investor to outperform by modulating the control $\gamma$ for a given control $h$. This can be conceptualized as a risk-sensitive zero sum stochastic differential game between the investor on one side and the market on the other and is formalized as follows\\
{\bf Problem (P1)}
Obtain $\hat{h} \in \mathcal{H}(T)$ and $\hat{\gamma} \in \Gamma(T)$ such that,
\begin{subequations}
\begin{align}
{J(f,x,\hat{h},\hat{\gamma};T)}
=\sup_{h \in {\mathcal{H}}(T)}\inf_{\gamma \in \Gamma(T)}{\frac{-2}{\theta}\log{E[(\frac{V^{h}_{T}}{L^{\gamma}_{T}})^{-\frac{\theta}{2}}]}}
=\inf_{\gamma \in \Gamma(T)}\sup_{h \in {\mathcal{H}}(T)}{\frac{-2}{\theta}\log{E[(\frac{V^{h}_{T}}{L^{\gamma}_{T}})^{-\frac{\theta}{2}}]}}
\end{align}
\end{subequations}
This can be construed as a game-theoretic version of the RSBAM problem.\\
{\bf Remark 2.1:}\\
~The problem set up ({\bf P1}) is an extension of Kuroda and Nagai  \cite{KN} and Davis and Lleo \cite{DL}. However the former does not consider the  benchmarked version i.e. the benchmark index is identically one in \cite{KN} while in Davis and Lleo \cite{DL} though  have a benchmarked portfolio criterion, they solve the one player optimization problem and not the two player saddle point problem. \\
\indent In light of the mathematical preliminaries just discussed, we formally elaborate the plan to solve the zero sum stochastic differential game ({\bf P1}). \\
\textit{Step 1}~ We reformulate the original objective criterion as a power utility function to an  exponential of an integral function.\\
 \textit{Step 2}~Define a new path   functional $I(f, x, {h},\gamma,t; T)$ (refer equation (\ref{3.9})) related to the exponential of the integral function.
 Define $\bar{u}(t,x)$ to be the upper-value function while $\underline{u}(t,x)$ be the lower-value function for the game associated with  $I$. Denote the game related to this objective functional as ({\bf P2}). \\
 \textit{Step 3}~Deduce the HJBI PDE corresponding to game ({\bf P2})( refer ({3.11}).\\
 \textit{Step 4}~Formulate the conditions that a candidate value function  should satisfy for the game with regards to objective function $I$ to have a value. This constitutes the verification lemma.\\
 \textit{Step 5}~Solve the HJBI PDE derived in step 3 while obtaining the expression for optimal controls. This optimal control pair will constitute a saddle point equilibrium for ({\bf P2}). The candidate value function satisfying all the conditions of the verification lemma is our desired value function for ({\bf P2}). \\
 \textit{Step 6}~Reverting back to the original problem {(\bf P1)}, show   using facts derived in Step 4, that  the game with objective criterion $J$ now has a value as well, and is in fact $u(0,x)$.\\
\indent In the next section we reformulate the objective criterion and formalize our game problem.
\section{Problem Reformulation}
\textit{Step 1}\\
We  will first transform the   utility optimization problem (\ref{2.5}) into optimizing the exponential-of-integral performance criterion.\\
{\bf Criterion under the expectation}\\
Our first aim is to write the objective criterion $J$ only in terms of the factor process. Towards that end we define the  function $g(x,h,\gamma,r;\theta)$ as follows:\\
\begin{eqnarray}\label{3.1}
g(x,h,\gamma,r;\theta)&=&\frac{1}{2}(\frac{\theta}{2}+1)h^{'}\Sigma\Sigma^{'}h-r-h^{'}(a+Ax-r1)+(\alpha+\beta x)-\frac{1}{2}{\frac{\theta}{2}}(h^{'}\Sigma\gamma+\gamma^{'}\Sigma^{'}h) \no\\ &+&\frac{1}{2}(\frac{\theta}{2}-1)\gamma^{'}\gamma
\end{eqnarray}
From (\ref{2.7})  and (\ref{3.1}) we therefore have,
\begin{eqnarray}\label{3.2}
 d\exp(\frac{-\theta}{2}F(t;h,\gamma))&=&\frac{\theta}{2}\bigg(g(X_{t},h_{t},\gamma_{t},r_{t};\theta)-(h_{t}^{'}\Sigma -\gamma^{'}_{t})\Sigma dW_{t}\bigg)-\frac{\theta^{2}}{8}(h^{'}_{t}\Sigma-\gamma^{'}_{t})\Sigma\Sigma^{'}(\Sigma^{'} h_{t}-\gamma_{t})dt\no\\
\end{eqnarray}
Thus we have,
\begin{eqnarray}\label{3.3}
\exp(\frac{-\theta}{2}F(t;h,\gamma))&=&f^{-\theta/2}\exp\{\frac{\theta}{2}\int_{0}^{t}{g(X_{s},h_{s},\gamma_{s},r;\theta)}ds\no \\ &-&\frac{\theta}{2}\int_{0}^{t}{(h^{'}_{s}\Sigma-\gamma^{'}_{s})}dW_{s}-\frac{1}{2}{(\frac{\theta}{2})}^{2}\int_{0}^{t}{(h^{'}_{s}\Sigma-\gamma^{'}_{s})(h^{'}_{s}\Sigma-\gamma^{'}_{s})^{'}}ds\}
\end{eqnarray}
where $V^{h}_{0}=v, L^{\gamma}_{0}=l$ and $f=\frac{V^{h}_{0}}{L^{\gamma}_{0}}=\frac{v}{l}$. \\
{\bf Change of measure}\\
 Let $\mathbb{P}^{h,\gamma}$ be the measure on ($\Omega,\mathcal{F}$) defined by,
\begin{eqnarray}\label{3.4}
\frac{d\mathbb{P}^{h,\gamma}}{d\mathbb{P}}|_{\mathcal{F}_{t}}&=& \bar{\mathcal{X}_{t}},
\end{eqnarray}
where $\bar{\mathcal{X}_{t}}$ is given by
\begin{eqnarray}\label{3.5}
\bar{\mathcal{X}_{t}}=\mathcal{E}(\frac{\theta}{2}\int_{0} {(h^{'}\Sigma-\gamma^{'})}dW)_{t}
\end{eqnarray}
and where $\mathcal{E}(\cdot)$ denotes the Doleans-Dade or martingale exponential. From the assumption made on the space of admissible controls $\mathcal{H}(T)$ and $\Gamma(T)$ it is clear that the Kazamaki condition $E[e^{\int_{0}^{t}{\theta{\frac{h^{'}_{s}\Sigma-\gamma^{'}_{s}}{2}}}dW_{s}}]<\infty$ $\forall t \in [0,T]$ is satisfied so that  $\mathbb{P}^{h,\gamma}$ to be a probability measure. i.e.
\begin{eqnarray}\label{3.6}
E[\mathcal{E}(\frac{\theta}{2}\int_{0} {(h^{'}\Sigma-\gamma^{'})}dW)_{T}]=1.
\end{eqnarray}
We note that,
\begin{eqnarray}\label{3.7}
W^{h,\gamma}_{t} \triangleq W_{t}+\frac{\theta}{2}\int_{0}^{t}{(h^{'}_{s}\Sigma-\gamma^{'}_{s})}ds,
\end{eqnarray}
 by Girsanov's formula, is a standard Brownian motion under $\mathbb{P}^{h,\gamma}$ and the factor process $X_{t}$ satisfies,
\begin{eqnarray}\label{3.8}
dX_{t}=(b+BX_{t}-\frac{\theta}{2}(\Sigma^{'} h_{t}-\gamma_{t}))^{'}dt+\Lambda dW^{h,\gamma}_{t}
\end{eqnarray}
\textit{Step 2}\\
{\bf The HJB equation}\\
Taking expectation w.r.t to the physical measure $\mathbb{P}$ and multiplying both sides of equation (\ref{3.3}) by $\frac{-2}{\theta}$  followed by the change of measure argument of (\ref{3.4}-\ref{3.5}) one considers the new path functional $I$ defined  as
\begin{eqnarray}\label{3.9}
I(f,x,h,\gamma,t,T)=\log{f}-\frac{2}{\theta}\log{E}^{h,\gamma}[\exp{\{\frac{\theta}{2}\int_{0}^{T-t}{g(X_{s},h_{s},\gamma_{s},r_{s+t};\theta)}ds\}}]
\end{eqnarray}
and then the upper-value function and lower-value function $\bar{u}$ and $\underline{u}$ respectively  for the   game corresponding to the new path functional $I$ are given by :
\begin{subequations}
\begin{align}
\bar{u}(t,x)=\sup_{h \in {\mathcal{H}}(T)}\inf_{\gamma \in \Gamma(T)}I(f,x,h,\gamma,t,T)\\
\underline{u}(t,x)=\inf_{\gamma \in \Gamma(T)}\sup_{h \in {\mathcal{H}}(T)}I(f,x,h,\gamma,t,T)\\
u(t,x)=\bar{u}(t,x)=\underline{u}(t,x)
\end{align}
\end{subequations}
If  a pair of controls satisfy (3.10c), then the game corresponding to the new path functional $I$ has the value $u$ and the pair of controls constitutes saddle point strategies for the game with regards to $I$.
Let the \textit{exponentially transformed} function $\tilde{I}$ be defined as $\tilde{I}=\exp(-\frac{\theta}{2}I)$ and   $\tilde{u}(t,x):=\exp(-\frac{\theta}{2}u(t,x))$. We now consider the problem of determining the saddle-point equilibrium for the game corresponding to the new path functional $\tilde{I}$. We call this problem (${\bf P2}$) and it is formally stated as follows:\\
{\bf Problem P2}
Obtain $\hat{h} \in \mathcal{H}(T)$ and $\hat{\gamma} \in \Gamma(T)$ such that,
\begin{eqnarray}\label{3.11}
\tilde{u}(t,x)&=&\inf_{h \in {\mathcal{H}}(T)}\sup_{\gamma \in \Gamma(T)}\tilde{I}(f,x,h,\gamma,t,T)\no\\
&=&\sup_{\gamma \in \Gamma(T)}\inf_{h \in {\mathcal{H}}(T)}\tilde{I}(f,x,h,\gamma,t,T)\no\\
&=& E^{\hat{h},\hat{\gamma}}[\exp\{\frac{\theta}{2}\int_{0}^{T-t}{g(X_{s},\hat{h}_{s},\hat{\gamma}_{s},r_{s+t};\theta)}ds\}f^{-\theta/2}]
\end{eqnarray}
We now provide a verification lemma for this game. Let us first define the process $Y^{h,\gamma}(t) $ by \\
$dY^{h,\gamma}(t)= \begin{pmatrix} dt \\ dX_t \end{pmatrix}=\begin{pmatrix} dt \\ 
(b+BX_{t}-\frac{\theta}{2}(h^{'}_{t}\Sigma-\gamma^{'}_{t}))dt+\Lambda dW^{h,\gamma}_{t} \end{pmatrix}$\\
Let $y \triangleq $($t,x$). The control process $h(t)=h(t,\omega)$ and $\gamma(t)=\gamma(t,\omega)$ for $\omega \in \Omega$ can be assumed to be Markovian.  Let $\mathcal{O}=(0,T) \times \mathbb{R}^{n}$.
Then the process $Y^{h,\gamma}(t)$  is a  Markov  process whose generator $\tilde{\mathcal{A}}^{h,\gamma}$ acting on a function $\tilde{u}(t,x) \in C^{2}_{0}([0,T] \times \mathbb{R}^{n})$ is given by,
\begin{eqnarray}\label{3.12}
\tilde{\mathcal{A}}^{h,\gamma}\tilde{{u}}(t,x) &=&\frac{\partial \tilde{u}(t,x)}{\partial t}+(b+Bx-\frac{\theta}{2}\Lambda(\Sigma^{'}h-\gamma))^{'}D{\tilde{u}}(t,x) +\frac{1}{2}tr(\Lambda\Lambda^{*}D^{2}{\tilde{u}}(t,x) )
\end{eqnarray}
in which $D\tilde{u}(t,x) \triangleq (\frac{\partial \tilde{u}(t,x)}{\partial x^{1}},...,\frac{\partial \tilde{u}(t,x)}{\partial x^{n}})^{'}$ and $D^{2}\tilde{u}(t,x)$ is the matrix defined by $D^{2}\tilde{u}(t,x)\triangleq [\frac{\partial^{2}\tilde{u}(t,x)}{\partial x^{i}x^{j}}],i,j=1,2,...,n.$\\
\textit{Step 3}\\
\indent By an application of the Feynman-Kac formula, it can be deduced from (\ref{3.11}) that the HJBI PDE for $\tilde{u}(t,x)$ is given by
\begin{eqnarray}\label{3.13}
{\bigg(\tilde{\mathcal{A}}^{\hat{h},\hat{\gamma}}+\frac{\theta}{2}g(x,\hat{h},\hat{\gamma},r;\theta)\bigg)}\tilde{u}(t,x)=0
\end{eqnarray}
Reversing the exponential transformation , dividing by $-(\theta/2)\tilde{u}(t,x)$, we can deduce from ({3.13}) that
 the HJBI PDE for $u(t,x)$ is given for $h \in \mathbb{R}^{m}$ and $\gamma \in \mathbb{R}^{(m+n)}$ by
\begin{eqnarray}\label{3.14}
 {\mathcal{A}}^{\hat{h},\hat{\gamma}}{u}(t,x)=0
\end{eqnarray}
 where the operator ${\mathcal{A}}^{{h},{\gamma}}$ is given by,
\begin{eqnarray}\label{3.15}
{\mathcal{A}^{h,\gamma}}{{u}}(t,x)&=&\frac{\partial {u}(t,x)}{\partial t} +(b+Bx-\frac{\theta}{2}\Lambda(\Sigma^{'}h-\gamma))^{'}Du(t,x)+\frac{1}{2}tr(\Lambda\Lambda^{'}D^{2}u(t,x))\no\\
&-&\frac{\theta}{4}(Du(t,x))^{'}\Lambda\Lambda^{'}Du(t,x)-g(x,h,\gamma,r;\theta)
\end{eqnarray}
In the next section we  provide a verification lemma for the game based on the  criterion function $I$. \\
\section{ Verification lemma for the game {\bf PII}}
\textit{Step 4}\\
We now provide a verification lemma related to the game ({\bf PII}).\\
{\bf Proposition 4.1.}~~\textit{Suppose $\tilde{w}$ $\in$ $\mathcal{C}^{1,2}({\mathcal{O}}) \cap \mathcal{C}(\bar{\mathcal{O}}) $ (is the space of twice differentiable functions on $\mathcal{O}$ with respect to $x$, once continuously differentiable on $\mathcal{O}$ with respect to $t$  and which are continuous on $\bar{\mathcal{O}}$ ). Suppose there exists a (Markov) control $\hat{h}(y),\hat{\gamma}(y)$ such that\\ 
1. $({\mathcal{\tilde{A}}}^{{h},\hat{\gamma}(y)}+\frac{\theta}{2}g(x,h,\hat{\gamma}(y),r;\theta))[(\tilde{w}(y))] \geq 0~ \forall~ h \in \mathbb{R}^{m}$;\\
2. $({\mathcal{\tilde{A}}}^{\hat{h}(y),{\gamma}}+\frac{\theta}{2}g(x,\hat{h}(y),\gamma,r;\theta))[(\tilde{w}(y))]  \leq 0~ \forall~\gamma \in  \mathbb{R}^{m+n}$;\\
3. $({\mathcal{\tilde{A}}}^{\hat{h}(y),\hat{\gamma}(y)}+\frac{\theta}{2}g(x,\hat{h}(y),\hat{\gamma}(y),r;\theta))[(\tilde{w}(y))]  = 0~ \forall~ y \in \mathcal{O}$;\\
4. $(\tilde{w}(T,X_{T}))={f}^{-\theta/2}$.\\
Define,
\begin{eqnarray}\label{4.1}
\tilde{Z}({s})&=&\tilde{Z}_{({s})}(h,\gamma)=\frac{\theta}{2}\bigg\{\int_{0}^{s}{g(X_{\tau},h_{\tau},\gamma_{\tau},r_{t+\tau};\theta)}d\tau\bigg\}
\end{eqnarray}
5. $E^{{h},{\gamma}}[\int_{0}^{T-t}{D\tilde{w}^{'}(t+s,X_{s})\Lambda}e^{\tilde{Z}_{s}}dW^{h,\gamma}_{s}]=0 ~ \forall~ h \in \mathbb{R}^{m},\forall~\gamma \in  \mathbb{R}^{m+n} $ \\
Now, define for each $y \in  \mathcal{O}$ and $h \in \mathcal{H}(T)$ and $\gamma\in \Gamma(T) $,
\begin{eqnarray*}
\tilde{I}(f,x,h,\gamma,t,T)&=&\exp(-\frac{\theta}{2}I(f,x,h,\gamma,t,T))\\
&=& E^{h,\gamma}[\exp\{\frac{\theta}{2}\int_{0}^{T-t}{g(X_{s},h_{s},\gamma_{s},r_{s+t};\theta)}ds\}f^{-\theta/2}],
\end{eqnarray*}
Then ($\hat{h}(y),\hat{\gamma}(y)$) is an optimal (Markov) control i.e.,
\begin{eqnarray*}
\tilde{w}(0,x)=\tilde{u}(0,x)=\tilde{I}({f,x,\hat{h},\gamma,0,T})& = & \inf_{h \in {\mathcal{H}}(T)}\{\sup_{\gamma \in \Gamma(T)}[\tilde{I}({f,x,{h},{\gamma},0,T})]\}  \\
&=& \sup_{\gamma \in \Gamma(T)}\{\inf_{h \in {\mathcal{H}}(T)}[\tilde{I}({f,x,{h},{\gamma},0,T})]\}  \\
&=& \sup_{\gamma \in \Gamma(T)}\tilde{I}({f,x,\hat{h},{\gamma},0,T})  \\
&=& \inf_{h \in {\mathcal{H}}(T)}\tilde{I}({f,x,{h},\hat{\gamma},0,T})=\tilde{I}({f,x,\hat{h},\hat{\gamma},0,T})
\end{eqnarray*}
}\\
{\bf Proof}~~
Apply  Ito's formula to $\tilde{w}(s,X_{s})e^{\tilde{Z}_{s}}$ to obtain
\begin{eqnarray}\label{4.2}
d(\tilde{w}(t+s,X_{s})e^{\tilde{Z}_{s}})&=&\bigg[e^{\tilde{Z}_{s}}({\mathcal{\tilde{A}}}^{{h},{\gamma}}+\frac{\theta}{2}g(X_{s},h_{s},\gamma_{s},r_{s+t};\theta))\bigg][(\tilde{w}(t+s,X_{s}))]ds+e^{\tilde{Z}_{s}}(D\tilde{w}(t+s,X_{s}))dW^{h,\gamma}_{s}\no\\ \no\\
\tilde{w}(T,X_{T-t})e^{\tilde{Z}_{T-t}}&=&\tilde{w}(t,x)+\int_{0}^{T-t}{((\tilde{\mathcal{A}}^{h,\gamma}+\frac{\theta}{2}g(X_{s},h_{s},\gamma_{s},r_{s+t};\theta))\tilde{w}(t+s,X_s))e^{\tilde{Z}_{s}}}ds\no\\
&+&\int_{0}^{T-t}{(D\tilde{w}^{'}(t+s,X_{s})\Lambda)e^{\tilde{Z}_{s}}}dW^{h,\gamma}_{s}
\end{eqnarray}
From condition(4) of statement of the Proposition, we have $\tilde{w}(T,X_T)=f^{-\theta/2}$. Taking expectation with respect to $\mathbb{P}^{h,\gamma}$ , setting $t=0$ and using conditions (1) and (5) of the Proposition  we get\\
\begin{eqnarray*}
E^{{h},\Gamma}[\tilde{w}(T,X_{T})e^{\tilde{Z}_{T}}]\geq \tilde{w}(0,x)
\end{eqnarray*}
Since this inequality is true for all $h \in {\mathcal{H}}(T)$ we have
\begin{eqnarray*}
\inf_{h \in {\mathcal{H}}(T)}E^{{h},\Gamma}[f^{-\theta/2}e^{\tilde{Z}_{T}}]\geq \tilde{w}(0,x)
\end{eqnarray*}
Hence we have,
\begin{eqnarray}\label{4.3}
\sup_{\gamma \in \Gamma(T)}\inf_{h \in {\mathcal{H}}(T)}E^{{h},{\gamma}}[f^{-\theta/2}e^{\tilde{Z}_{T}}]\geq \inf_{h \in {\mathcal{H}}(T)}E^{{h},\Gamma}[f^{-\theta/2}e^{\tilde{Z}_{T}}]\geq  \tilde{w}(0,x)
\end{eqnarray}
Similarly, setting $t=0$ we get, using condition (2) of the Proposition, we get  the following lower bound,\\
\begin{eqnarray*}
E^{\hat{h},{\gamma}}[\tilde{w}(T,X_{T})e^{\tilde{Z}_{T}}]\leq \tilde{w}(0,x)
\end{eqnarray*}
Since this inequality is true for all $\gamma \in \Gamma(T)$ we have
\begin{eqnarray*}
\sup_{\gamma \in \Gamma(T)}E^{\hat{h},{\gamma}}[f^{-\theta/2}e^{\tilde{Z}_{T}}]\leq  \tilde{w}(0,x)
\end{eqnarray*}
Hence we have,
\begin{eqnarray}\label{4.4}
\inf_{h \in {\mathcal{H}}(T)}\sup_{\gamma \in \Gamma(T)}E^{{h},{\gamma}}[f^{-\theta/2}e^{\tilde{Z}_{T}}]\leq
\sup_{\gamma \in \Gamma(T)}E^{\hat{h},{\gamma}}[f^{-\theta/2}e^{\tilde{Z}_{T}}]\leq \tilde{w}(0,x)
\end{eqnarray}
Also , setting $t=0$ and using condition (3) of the Proposition and using the definition of  $\tilde{u}$ in (\ref{3.11}) we get,\\
\begin{eqnarray}\label{4.5}
E^{\hat{h},\hat{\gamma}}[\tilde{w}(T,X_{T})e^{\tilde{Z}_{T}}]&=& \tilde{w}(0,x)\no\\
&=&E^{\hat{h},\hat{\gamma}}[\exp\{\frac{\theta}{2}\int_{0}^{T}{g(X_{s},\hat{h}_{s},\hat{\gamma}_{s},r_{s+t};\theta)}ds\}f^{-\theta/2}]
\end{eqnarray}
It is automaticaly true that
\begin{eqnarray}\label{4.6}
\sup_{\gamma \in \Gamma(T)}\inf_{h\in {\mathcal{H}}(T)}E^{{h},{\gamma}}[f^{-\theta/2}e^{\tilde{Z}_{T}}]\leq \inf_{h\in {\mathcal{H}}(T)}\sup_{\gamma\in \Gamma(T)}E^{{h},{\gamma}}[f^{-\theta/2}e^{\tilde{Z}_{T}}].
\end{eqnarray}
Conversely, from (\ref{4.3}), (\ref{4.4}) and (\ref{4.5}) we have,
\begin{eqnarray}\label{4.7}
\inf_{h \in {\mathcal{H}}(T)}\sup_{\gamma \in \Gamma(T)}E^{{h},{\gamma}}[f^{-\theta/2}e^{\tilde{Z}_{T}}]&\leq & \tilde{w}(0,x) \leq  \sup_{\gamma \in \Gamma(T)}\inf_{h\in {\mathcal{H}}(T)}E^{{h},{\gamma}}[f^{-\theta/2}e^{\tilde{Z}_{T}}]
\end{eqnarray}
Hence from (\ref{4.6}) and (\ref{4.7}) we have,
\begin{eqnarray}\label{4.8}
\sup_{\gamma \in \Gamma(T)}\inf_{h\in {\mathcal{H}}(T)}E^{{h},{\gamma}}[f^{-\theta/2}e^{\tilde{Z}_{T}}]&= & \inf_{h \in {\mathcal{H}}(T)}\sup_{\gamma \in \Gamma(T)}E^{{h},{\gamma}}[f^{-\theta/2}e^{\tilde{Z}_{T}}]\no\\
&=& \tilde{w}(0,x) = E^{\hat{h},\hat{\gamma}}[f^{-\theta/2}e^{\tilde{Z}_{T}}]
\end{eqnarray}
$\qed$\\\\
{\bf Corollary 4.2}~~\textit{Admissible(optimal) strategies for the exponentially transformed problem  given by (\ref{3.11}) are also admissible(optimal) for the problem (3.10c). Formally,
\begin{eqnarray*}
u(0,x)& = & \sup_{h \in {\mathcal{H}}(T)}\{\inf_{\gamma \in \Gamma(T)}[{I}({f,x,{h},{\gamma},0,T})]\}  \\
&=& \inf_{\gamma \in \Gamma(T)}\{\sup_{h \in {\mathcal{H}}(T)}[{I}({f,x,{h},{\gamma},0,T})]\}  \\
&=& \inf_{\gamma \in \Gamma(T)}{I}({f,x,\hat{h},{\gamma},0,T})  \\
&=& \sup_{h \in {\mathcal{H}}(T)}{I}({f,x,{h},\hat{\gamma},0,T})={I}({f,x,\hat{h},\hat{\gamma},0,T})
\end{eqnarray*}
}
{\bf Proof}~~ The value function $u$ and $\tilde{u}$ are related through the strictly monotone continuous transformation $\tilde{u}(t,x)=\exp(-\frac{\theta}{2}u(t,x))$. Thus admissible (optimal) strategies for the exponentially transformed problem are also admissible(optimal) for the  problem (3.10c). $\qed$
\section{Solving the risk-sensitive zero sum stochastic differential game}
\textit{Step 5}\\
\indent We seek to find the  value function $u$ for the game defined in (\ref{3.12}). We guess a solution  assuming that it belongs to the class $C^{1,2}((0,T) \times \mathbb{R}^{n})$  and show that the guess satisfies all the conditions of our verification lemma given by Proposition 4.1.  Conditions (1)-(4) of the verification lemma can be written in a compact form as
\begin{eqnarray}\label{5.1}
\sup_{h \in {\mathcal{H}}(T)}\inf_{\gamma \in \Gamma(T)}{\mathcal{A}}^{{h},{\gamma}}{u}(t,x)=0;~~~ {u}(T,x)=\log{f}
\end{eqnarray}
Motivated by the results in Kuroda and Nagai \cite{KN}, we will  look for a $u$ given by $u(t,x)=\frac{1}{2}x^{'}Q_{t}x+ q^{'}_{t}x + k_{t}$ where $Q$ is an $n \times n$ symmetric matrix, $q \in \mathbb{R}^{n}$ and $k$ is a scalar. Substituting this form in (\ref{3.15}) we get
\begin{eqnarray}\label{5.2}
\mathcal{A}^{h,\gamma}u(t,x)&=&\frac{1}{2}x^{'}\frac{d Q_{t}}{dt}x+{\frac{dq_{t}}{dt}}^{'}x+\frac{dk_{t}}{dt}+\bigg(b+Bx-\frac{\theta}{2}\Lambda(\Sigma^{'}h_{t}-\gamma(t))\bigg)^{'}(Q_{t}x+q_{t}) \no\\
&+&\frac{1}{2}(\Lambda\Lambda^{'}Q_{t}Q^{'}_{t}\Lambda^{'}\Lambda)-\frac{\theta}{4}(Q_{t}x+k_{t})^{'}\Lambda\Lambda^{'}(Q_{t}x+k_{t})\no\\
&-&\frac{1}{2}(\frac{\theta}{2}+1)h^{'}_{t}\Sigma\Sigma^{'}h_{t}+r_{t}-(\alpha_{t}+\beta x)+h^{'}_{t}(a+Ax-r_{t}1)+\frac{1}{2}{\frac{\theta}{2}}(h^{'}_{t}\Sigma\gamma+\gamma^{'}\Sigma^{'}h_{t})\no\\
&-&\frac{1}{2}(\frac{\theta}{2}-1)\gamma^{'}_{t}\gamma_{t}
\end{eqnarray}
{\bf Remark 5.1}~~\textit{ Since the game considered is for the risk-averse investor $ \theta>0$. Moreover based in the expression for $\hat{\gamma}$ in (\ref{5.5}), $\theta \neq 2$. This leaves for two possibilities: $\theta \in (0,2)$ or $ \theta \in (2, \infty)$. For the optimal strategies ($\hat{h},\hat{\gamma}$) to be a saddle-point equilibrium for the game, we would desire that the equation with the quadratic term in $h$ be negative definite while the quadratic term in $\gamma$ be positive definite. In fact for the choice  $\theta>0$, the quadratic term in $h$ desirably is negative definite  while for $\theta<2$, the quadratic term in $\gamma$ is positive definite . Hence for our case the valid range of $\theta$ is between 0 and 2 and excludes the other two possibilities for the range of $\theta$.}\\
We now  solve the first order condition for $\hat{\gamma}$ to minimize ${\mathcal{A}}^{\hat{h},{\gamma}}{u}(t,x)$ over all $\gamma \in \mathbb{R}^{n+m}$:
\begin{eqnarray}\label{5.3}
({2-\theta})\hat{\gamma_{t}}-\theta(\Sigma^{'}\hat{h}_{t}-\hat{\gamma}^{'})Du(t,x)=0
\end{eqnarray}
The first order condition for $\hat{h}$ that maximizes ${\mathcal{A}}^{{h},\hat{\gamma}(y)}\tilde{u}(t,x)$ over all $h \in \mathbb{R}^{m}$ in terms of $u(t,x)$ is,
\begin{eqnarray}\label{5.4}
\hat{h}_{t}=\frac{2}{(\theta+2)}(\Sigma\Sigma^{'})^{-1}[d_{t}+\frac{\theta}{2}\Sigma\hat{\gamma}_{t}-\frac{\theta}{2}\Sigma\Lambda^{'}D u(t,x)]
\end{eqnarray}
Substituting back $\hat{h}$ obtained in (\ref{5.4}) into (\ref{5.3}) we get
\begin{eqnarray}\label{5.5}
\hat{\gamma}_{t}&=&\frac{\theta}{2-\theta}[\Sigma^{'}\hat{h}_{t}-\Lambda^{'}Du(t,x)]
\end{eqnarray}
The optimal control $\hat{h}_t$ is a global maximum while $\hat{\gamma}_t$ is a global minimum for $t \leq [0,T]$.
We substitute $\hat{h}$ from (\ref{5.4}) and $\hat{\gamma}$ from (\ref{5.5}) in (\ref{5.1}) to obtain
\begin{eqnarray}\label{5.6}
\mathcal{A}^{\hat{h},\hat{\gamma}}u(t,x)=0; ~~~~{u}(T,x)=\log{f}
\end{eqnarray}
We then group all the resulting quadratic terms in $x$, linear terms in $x$ and constants together to  conclude that the choice of $u(t,x)=\frac{1}{2}x^{'}Q_{t}x+ q^{'}_{t}x + k_{t}$ is indeed the solution to the HJBI PDE (\ref{5.1})
provided that $Q$, $q$ and $k$ satisfy the following system of differential equations:\\
$\bullet$ a matrix Ricatti equation related to the coefficient of the quadratic term and used to determine the symmetric non-negative matrix $Q_{t}$, given as
\begin{eqnarray}\label{5.7}
\frac{d Q_{t}}{dt}&=&Q_{t}K_{0}Q_{t}+K_{1}^{'}Q_{t}+Q_{t}K_{1}+2\frac{2-\theta}{{(2-\theta^2)}^2}A^{'}{(\Sigma\Sigma^{-1})}^{-1}A=0~~~0 \leq t \leq T,\no\\~~
Q_{T}&=&0
\end{eqnarray}
{where}\\
$K_{0}=\frac{-\theta^2}{2(2-\theta)} \Lambda \Lambda^{'}+\frac{2\theta^2}{(2-\theta){(2-\theta^2)}^2}\Lambda \Sigma^{'}{(\Sigma\Sigma^{'})}^{-1}\Sigma\Lambda^{'}$\\
$K_{1}= B-\frac{2\theta}{{(2-\theta^2)}^2}A^{'}{(\Sigma\Sigma^{'})}^{-1}\Sigma\Lambda^{'}$\\
$\bullet$~~ The following linear ordinary differential equation satisfied by the $n$ element column vector $q(t)$
\begin{eqnarray}\label{5.8}
\frac{dq_{t}}{dt}&+&(K_{1}^{'}+Q_{t}K_{0})q_{t}+Q^{'}_{t}b+(a-r(t)1)^{'}{(\Sigma\Sigma^{'})}^{-1}[\frac{-2\theta}{{(2-\theta^2)}^2}\Sigma\Lambda^{'}Q(t)+\frac{(2-\theta)}{(2-\theta^2)^2}A]\no\\
&-&\beta_t\no \\
&&q_{T}=0
\end{eqnarray}
$\bullet$~~ The following linear ordinary differential equation satisfied by the constant $k_{t}$
\begin{eqnarray}\label{5.9}
&&\frac{dk_{t}}{dt}+\frac{1}{2}tr(\Lambda\Lambda^{'}Q_{t})+r_t-\alpha_t-\frac{2\theta}{(2-\theta^2)^2}(a-r(t)1)^{'}{(\Sigma\Sigma^{'})}^{-1}\Sigma\Lambda^{'}q(t)\no\\
&+&\frac{2-\theta}{(2-\theta^2)^2}(a-r(t)1)^{-1}{(\Sigma\Sigma^{'})}^{-1}(a-r(t)1)+\frac{\theta^2}{(2-\theta)(2-{\theta}^2)^2}q^{'}(t)\Lambda\Sigma^{'}{(\Sigma\Sigma^{'})}^{-1}\Sigma\Lambda^{'}q(t)\no\\
&-& \frac{\theta^2}{4(2-\theta)}q^{'}(t)\Lambda\Lambda^{'}q(t)\no\\
k_{T}&=&\log{f}
\end{eqnarray}
Condition 4 of Proposition 4.1 in terms of $u$ imposes the terminal condition in (\ref{5.9}). \\
\indent If $K_{0}$ is positive definite then a unique solution to the Riccati equation (\ref{5.7}), $Q_{t}$ , exists  for all $t \leq T$. This property of  positive definiteness follows from interpretation of the solution $Q_{t}$ as the covariance matrix of observations from a Kalman filter used to estimate the state of a dynamical system (see Theorem 4.4.1 in Davis \cite{Davis}) for details. The uniqueness  property of $Q_{t}$ follows from the standard existence-uniqueness theorem for first order differential equations (see Proposition 4.4.2 in Davis \cite{Davis}).   \\
It remains to be seen if $\tilde{u}=\exp(-\frac{\theta}{2}u)$ for the choice of $u$ satisfies condition (5) of Proposition 4.1. \\
 {\bf Proposition 5.2}~~\textit{$E^{{h},{\gamma}}[\int_{0}^{T-t}{e^{\tilde{Z}_{s}}(D\tilde{u}^{'}(t+s,X_{s})\Lambda)}dW^{h,\gamma}_{s}]=0$.}\\
 {\bf Proof}~~From the definition of $\tilde{u}$ in (\ref{3.11}), for any optimal control belonging to $\Gamma(T)$, the strategy $\hat{h} \equiv 0$ is sub-optimal, and hence will provide an upper bound on $\tilde{u}$. Further for the zero-benchmark case namely, $\hat{\gamma} \equiv 0$, we would obtain now an upper bound on $\tilde{u}$
 \begin{eqnarray*}
\tilde{u}(t,x)&=&\inf_{h \in {\mathcal{H}}(T)}E^{h,\hat{\gamma}}[\exp\{\frac{\theta}{2}\int_{0}^{T-t}{g(X_{s},h_{s},\hat{\gamma}_{s},r_{s+t};\theta)}ds\}f^{-\theta/2}]\\
&\leq & E^{0,\hat{\gamma}}[\exp\{\frac{\theta}{2}\int_{0}^{T-t}{g(X_{s},0,\hat{\gamma}_{s},r_{s+t};\theta)}ds\}f^{-\theta/2}]\\
\therefore \tilde{u}(t,x)&\leq & E^{0,0}[\exp\{\frac{\theta}{2}\int_{0}^{T-t}{g(X_{s},0,0,r_{s+t};\theta)}ds\}f^{-\theta/2}]\\
&= & \exp(-\frac{\theta}{2}\int_{0}^{T-t}{r_{s+t}}ds)f^{-\theta/2}
 \end{eqnarray*}
Now $Q$ and $q$ are solutions to the system of o.d.e, and  hence are  integrals of  bounded functions  . Hence $Q$ and $q$  are continuous functions of time $t \in [0,T]$ and hence bounded on $[0,T]$. The matrix $\Lambda$ is a known constant. From standard existence-uniqueness result of stochastic differential equation (refer Oksendal (\cite{Oks})) we have $X \in L^{2}(\Omega,\mathcal{F},\mathbb{P}^{h,\gamma})$. Hence from the upper bound on $\tilde{u}$ , Remark 3.1 and the fact that $Du(t,X_{t})=Q_{t}X_{t}+q_{t}$ is in $L^{2}(\Omega,\mathcal{F},\mathbb{P}^{h,\gamma})$, we have that $E^{h,\gamma}([D\tilde{u}~\Lambda e^{\tilde{Z}}, D\tilde{u}~\Lambda e^{\tilde{Z}}]_{t})< \infty$ $\forall t \in [0,T]$. Hence we have $E^{{h},{\gamma}}[\int_{0}^{T-t}{D\tilde{u}^{'}(t+s,X_{s})\Lambda}e^{\tilde{Z}_{s}}dW^{h,\gamma}_{s}]=0$. $\qed$.\\
It is clear that our guess for $\tilde{u}=\exp(-\frac{\theta}{2}u)$ satisfies conditions (1)-(5) of Proposition 4.1. Hence our choice of $\tilde{u}$  indeed is the  value of the game ({\bf P2}) and controls $\hat{h},\hat{\gamma}$ are the saddle point equilibrium of this game.
\\
{\bf Lemma 5.2}~\textit{For the choice of space of controls  $\mathcal{H}(T)$ and $\Gamma(T)$, we have
\begin{eqnarray}\label{5.10}
E[\mathcal{E}\bigg(-\frac{\theta}{2}\int_{0}{[(Q_{t}X_{t}+q_{t})\Lambda+({h}^{'}_{t}\Sigma-\gamma^{'}_{t})]dW_{t}}\bigg)_{T}]=1
\end{eqnarray}}
{\bf Proof:} From the Kazamaki condition, refer (Oksendal \cite{Oks}), ({5.10}) holds if \\ $E[\exp(\int_{0}^{t}{\theta(\frac{(Q_{s}X_{s}+q_{s})\Lambda+({h}^{'}_{s}\Sigma-\gamma^{'}_{s})}{2})}dW_{s})]<\infty$ $\forall~t \in [0,T]$. Hence by application of Cauchy-Schwartz inequality we have,
\begin{eqnarray*}
E[\exp(\int_{0}^{t}{\theta(\frac{(Q_{s}X_{s}+q_{s})\Lambda+({h}^{'}_{s}\Sigma-\gamma^{'}_{s})}{2})}dW_{s})] &\leq&(E[e^{\int_{0}^{t}{\theta{(Q_{s}X_{s}+q_{s})\Lambda}}dW_{s}}])^{1/2}{(E[e^{\int_{0}^{t}{\theta{{({h}^{'}_{s}\Sigma-\gamma^{'}_{s})}}}dW_{s}}])}^{1/2}\no\\
\end{eqnarray*}
However for $E[e^{\int_{0}^{t}{\theta{(Q_{s}X_{s}+q_{s})\Lambda}}dW_{s}}]<\infty$ to hold , it is enough to show that the Novikov condition given by $E[e^{\int_{0}^{T}{\theta^{2}{(Q_{s}X_{s}+q_{s})\Lambda\Lambda^{'}(Q_{s}X_{s}+q_{s})}}d{s}}]<\infty$ hold; refer (Oksendal \cite{Oks}). Since $X$ is Gaussian process and $Q_{t}$ and $q_{t}$ are deterministic, $(Q_{t}X_{t}+q_{t})\Lambda$ is Gaussian and hence by completion of squares argument detailed in Theorem 5.3 below we have $E[e^{\int_{0}^{T}{\theta^{2}{(Q_{s}X_{s}+q_{s})\Lambda\Lambda^{'}(Q_{s}X_{s}+q_{s})}}d{s}}]<\infty$ holds and hence $E[e^{\int_{0}^{t}{\theta{(Q_{s}X_{s}+q_{s})\Lambda}}dW_{s}}]<\infty$  $\forall t \in [0,T]$ is validated. ${(E[e^{\int_{0}^{t}{\theta{{({h}^{'}_{s}\Sigma-\gamma^{'}_{s})}}}dW_{s}}])}^{1/2}< \infty$ is  validated from similar application of Cauchy-Schwartz inequality followed by the assumption made earlier in the definition of the space of controls $\mathcal{H}(T)$ and $\Gamma(T)$. Thus the Kazamaki condition holds and the conclusion follows. $\qed$\\
{\bf Theorem 5.3}~~\textit{If there exist a solution $Q$ to (\ref{5.7}), then the  strategies $(\hat{h},\hat{\gamma}) $ defined by
\begin{eqnarray}\label{5.10}
\hat{h}_{t}=\frac{2}{(\theta+2)}(\Sigma\Sigma^{'})^{-1}[d_{t}+\frac{\theta}{2}\Sigma\gamma_{t}-\frac{\theta}{2}\Sigma\Lambda^{'}(Q_{t}X_{t}+q_{t})]
\end{eqnarray}
\begin{eqnarray}\label{5.11}
\hat{\gamma}_{t}=\frac{\theta}{2-\theta}[\Sigma^{'}\hat{h}_{{t}}-\Lambda^{'}(Q_{t}X_{t}+q_{t})]
\end{eqnarray}
where $q$ is a solution of (\ref{5.8}) are admissible i.e. $ h \in {\mathcal{H}}(T)$ and $\gamma \in \Gamma(T)$ and are optimal for the finite horizon game problem ({\bf P1}), namely,
\begin{eqnarray*}
u(0,x)&=& \sup_{h \in {\mathcal{H}}(T)}\inf_{\gamma \in \Gamma(T)}{J}(f,x,h,\gamma,T;\theta)\\
&=&\inf_{\gamma \in \Gamma(T)}\sup_{h \in {\mathcal{H}}(T)}{J}(f,x,h,\gamma,T;\theta)\\
&=&\inf_{\gamma \in \Gamma(T)}{J}(f,x,\hat{h},\gamma,T;\theta)\\
&=& \sup_{h \in {\mathcal{H}}(T)}{J}(f,x,h,\hat{\gamma},T;\theta)\\
&=&{J}(f,x,\hat{h},\hat{\gamma},T;\theta)\\
&=&\frac{1}{2}x^{'}Q_0x+q^{'}_0x+k_0
 \end{eqnarray*}}
{\bf Proof}~~  The controls derived in section 5, $(\hat{h},\hat{\gamma})$ forms the saddle point equilibrium for the ({\bf P2}) game . We aim to show  that these controls are in fact admissible and  optimal for the problem ({\bf P1}) as well.\\
\textit{ Proof of admissibility}
From the expression for $\hat{h}$ and $\hat{\gamma}$ in (\ref{5.10}) and (\ref{5.11}) respectively we note that $ -\frac{\theta}{2}\bigg((Q_{t}X_{t}+q_{t})\Lambda+(\hat{h}^{'}_{t}\Sigma-\hat{\gamma}^{'}_{t})\bigg)$
can be written linearly in $X_t$ as $ X^{'}_{t}v^{1}_{t}+v^{2}_{t}$ where, constants $v^{1}_{t}$ and $v^{2}_{t}$ are given by,
\begin{eqnarray*}
v^{1}_{t}&=&-\frac{\theta}{2}Q^{'}(t)\Lambda+\frac{\theta(\theta-1)}{(2-\theta^2)}A^{'}(\Sigma\Sigma^{'})^{-1}\Sigma \Lambda^{'}+\frac{\theta(\theta-1)}{2-\theta^2}Q^{'}(t)\Lambda\Sigma^{'}(\Sigma\Sigma^{'})^{-1}(a-r1)\\
&-&\frac{2\theta^{2}(\theta-1)}{(2-\theta)(2-\theta^{2})}Q^{'}(t)\Lambda\Sigma^{'}(\Sigma\Sigma^{'})^{-1}\Sigma\Lambda^{'}q(t)-\frac{\theta^{2}}{(2-\theta)}Q^{'}(t)\Lambda\Lambda^{'}q(t).\\
v^{2}_t&=&-\frac{\theta}{2}q^{'}(t)\Lambda+\frac{\theta(\theta-1)}{(2-\theta^{2})}(a-r1)^{'}(\Sigma\Sigma^{'})^{-1}\Sigma\Lambda^{'}q(t)\\
&-&\frac{\theta^2(\theta-1)}{(2-\theta)(2-\theta^2)}q^{'}(t)\Lambda\Sigma^{'}(\Sigma\Sigma^{'})^{-1}\Sigma\Lambda^{'}q(t)-\frac{\theta^{2}}{(2-\theta)}q^{'}(t)\Lambda\Lambda^{'}q(t)
\end{eqnarray*}
Since $X$ satisfies the SDE , $dX_t=(b+BX_t)dt+\Lambda dW_t$, so~ $ E|X_t| \leq E|X(0)|+|b|T+|B|\int_{0}^{t}{E{|X_s|}ds}$. By Gronwall's inequality, therefore $E|X_t|\leq (E|X(0)|+|b|T)\exp(|B|t)$ and $Cov(X_t)=\Lambda^{'}\Lambda t$. Let $\phi(t)\triangleq v^{1}_{t}X_t+v^{2}_{t}$. We now  explicitly calculate $E[e^{\delta|\phi_{t}|^{2}}]$ for some $\delta>0$ since from Remark 2 in Lemma 2, of section 12 (Gihman and Skorokhod \cite{GS}) would imply that the Novikov's condition holds true. Let $R_{t}=e^{-Bt}X_t+e^{-bt}$.  Hence $dR_t=e^{-Bt}\Lambda dW_t$. Therefore $R_t$ is a Gaussian process and hence $\phi_t$ is Gaussian process with drift.
Also $\mu_t=E[|\phi_t|]\leq \sup_{0\leq t \leq T}|v^{1}_{t}|(E|X_0|+|b|T)\exp(|B|t)+\sup_{0\leq t \leq T}|v^{2}_{t}|$ and $\tilde{\Sigma}_t=Cov(\phi_t) \leq {v^{1}}^{'}_{t}\Lambda^{'}\Lambda v^{1}_{t}$. Thus mean $\mu_t$ and co-variance $\tilde{\Sigma}_t$ are bounded above by $t$.
We use the following completion of squares argument: $\frac{1}{2}z^{'}Az+b^{'}z+c=\frac{1}{2}(z+A^{-1}b)^{'}A(z+A^{-1}b)+c-\frac{1}{2}b^{'}A^{-1}b$ .
\begin{eqnarray*}
E[e^{\delta|\phi_{t}|^{2}}]&=&\int_{\mathbb{R}^{n}}{\frac{1}{{2\pi}^{n/2}|\tilde{\Sigma}_{t}\tilde{\Sigma}^{'}_{t}|^{1/2}}e^{\delta|\phi|^{2}_t}e^{-\frac{1}{2}(\phi-\mu_t)^{'}(\tilde{\Sigma}_{t}\tilde{\Sigma}^{'}_{t})^{-1}(\phi-\mu_{t})dx^{1}dx^{2}...dx^{n}}}\\
&=&\frac{1}{{2\pi}^{n/2}|\tilde{\Sigma}_{t}\tilde{\Sigma}^{'}_{t}|^{1/2}}\int_{\mathbb{R}^{n}}{e^{\frac{-{\phi^{'}{(-2\delta I+(\tilde{\Sigma}_{t}\tilde{\Sigma}^{'}_{t})^{-1})}^{-1}\phi+2\mu^{'}(t)(\tilde{\Sigma}_{t}\tilde{\Sigma}^{'}_{t})^{-1}\phi-\mu^{'}_t(\tilde{\Sigma}_{t}\tilde{\Sigma}^{'}_{t})^{-1}\mu_t}}{2}}}dx^{1}....dx^{n}\\
&=&\frac{|{(\tilde{\Sigma}^{'}_{t}\tilde{\Sigma}_{t})}|^{-1/2}}{{|(-2\delta I+(\tilde{\Sigma}_{t}\tilde{\Sigma}^{'}_{t})^{-1})^{-1}|}^{-1/2}}\times\\
&&e^{\frac{{-\mu^{'}_t(\tilde{\Sigma}_{t}\tilde{\Sigma}^{'}_{t})^{-1}\mu_t}+ 4\mu^{'}_t(\tilde{\Sigma}_{t}\tilde{\Sigma}^{'}_{t})^{-1}{(-2\delta I+(\tilde{\Sigma}_t\tilde{\Sigma}^{'}_t)^{-1})}^{-1}(\tilde{\Sigma}_t\tilde{\Sigma}^{'}_t)^{-1}\mu_{t}}{2}}
\end{eqnarray*}
Matrix ${(\tilde{\Sigma}_t\tilde{\Sigma}_t)}^{-1}$ is symmetric positive definite with lowest eigenvalue say $\lambda_{min}$. Then it is easy to show that for $\delta<\frac{\lambda_{min}}{2}$, matrix $(-2\delta I+(\tilde{\Sigma}_{t}\tilde{\Sigma}^{'}_{t})^{-1})^{-1}$ is positive definite . Along with the derived fact that  $\mu_t$ and $\tilde{\Sigma}_t$ is bounded above by $t\leq T$ , hence there exists some constant $C$ such that $ E[e^{\delta|\phi_{t}|^{2}}]\leq  C $. Hence the optimal controls ${\hat{h}},\hat{\gamma}$ belong to their respective admissible class viz. $\mathcal{H}(T)$ and $\Gamma(T)$ respectively. $\qed$\\
\textit{ Proof of optimality}~~ Define,
\begin{eqnarray}\label{5.12}
Z_{s}&=&Z_{s}(h,\gamma)=\frac{\theta}{2}\bigg\{\int_{0}^{s}{g(X_{\tau},h_{\tau},\gamma_{\tau},r_{t+\tau};\theta)}d\tau-{(h^{'}_{\tau}\Sigma-\gamma^{'}_{\tau})}dW_{\tau}\no\\
&-&\frac{\theta}{4}{(h^{'}_{\tau}\Sigma-\gamma^{'}_{\tau})}^{'}{(h^{'}_{\tau}\Sigma-\gamma^{'}_{\tau})}d\tau\bigg\}
\end{eqnarray}
Also define, $\chi(t,x)=-\frac{\theta}{2}(u(t,x)-\log{f})$ and $Lu(t,x)=\frac{1}{2}tr(\Lambda\Lambda^{'}D^{2}u(t,x))+(b+Bx)^{'}Du(t,x)$\\
Hence, we have\\
\begin{eqnarray*}
d\chi(t+s,X_s)&=&-\frac{\theta}{2}(\frac{\partial u}{\partial t}+Lu)(t+s,X_s)ds-\frac{\theta}{2}Du(t+s,X_s)^{'}\Lambda dW_s
\end{eqnarray*}
Hence,
\begin{eqnarray*}
\frac{d\exp\{\chi(t+s,X_s)\}}{\exp\{\chi(t+s,X_{s})\}}&=& -\frac{\theta}{2}(\frac{\partial u}{\partial t}(t,x)+{L}u)(t+s,X_s)-\frac{\theta}{2}Du(t+s,X_s)^{'}\Lambda dW_s\\
&+&\frac{\theta^{2}}{8}Du^{'}\Lambda\Lambda^{'}Du(t+s,X_{s})ds
\end{eqnarray*}
and so,
\begin{eqnarray*}
\frac{d\exp\{\chi(t+s,X_s)\}\exp\{Z(s)\}}{\exp\{\chi(t+s,X_s)\}\exp\{Z(s)\}}&=& -\frac{\theta}{2}(\frac{\partial u}{\partial t}(t,x)+{L}u)(t+s,X_s)
-\frac{\theta}{2}{Du(t+s,X_s)}^{'}\Lambda dW_s\\
&+&\frac{\theta^{2}}{8}Du^{'}\Lambda\Lambda^{'}Du(t+s,X_s)ds +\frac{\theta}{2}{g(X_s,h_s,\gamma_s,r_s+t;\theta)}ds\\
&-&\frac{\theta}{2}(h^{'}(s)\Sigma-\gamma^{'}(s))dW_s
+\frac{\theta^{2}}{4}(h^{'}(s)\Sigma-\gamma^{'}(s))\Lambda^{'}Du(t+s,X_s)ds
\end{eqnarray*}
Hence from (\ref{3.15}), we have,
\begin{eqnarray}\label{5.13}
&&\exp\{\chi(T,X(T-t))+Z(T-t)\}=\exp(\chi(t,x))\exp \bigg[\int_{0}^{T-t}{-\frac{\theta}{2}(\mathcal{A}^{h,\gamma} u(t+s,X_s))}ds\no\\
&-&\int_{0}^{T-t}{\frac{\theta}{2}[Du(t+s,X_s)^{'}\Lambda+(h^{'}_t\Sigma-\gamma^{'}_t)]}dW_t\no\\
&-&\int_{0}^{T-t}{\frac{{\theta}^{2}}{8}{[Du(t+s,X_s)^{'}+(h^{'}_t\Sigma-\gamma^{'}_t)][Du(t+s,X_s)^{'}+(h^{'}_t\Sigma-\gamma^{'}_t)]^{'}}ds}\bigg]\no\\
\end{eqnarray}
We have shown that $u$ satisfies  conditions (1)-(5) of Proposition 4.1 Hence from condition(4) of Proposition 4.1,  we have $\chi(T,x)=0$. Now  setting $t=0$ and taking condition (1) of  Proposition 4.1 into account for $\gamma=\hat{\gamma}$, and for any $h \in \hat {\mathcal{H}}(T)$ we see from (\ref{5.13}) that \\
\begin{eqnarray*}
{(\frac{V^{h}_{T}}{L^{\gamma}_{T}})}^{-\theta/2}\geq e^{-\frac{\theta}{2}u(0,x)}\exp\bigg[-\int_{0}^{T}{\frac{\theta}{2}[Du(s,X_s)^{'}\Lambda+(h^{'}_{s}\Sigma-\Gamma^{'}_s)]}dW_s\no\\
-\int_{0}^{T}{\frac{{\theta}^{2}}{8}{[Du(s,X_s)^{'}+(h^{'}_s \Sigma-\gamma^{'}_s)][Du(s,X_s)^{'}+(h^{'}_s\Sigma-\Gamma^{'}_s)]^{'}}ds}\bigg]
\end{eqnarray*}
Now by taking expectations w.r.t to the physical probability measure $\mathbb{P}$ on both sides of above equation  and using Lemma 5.2, we obtain\\
\begin{eqnarray*}
J(f,x,h,\gamma,T) \leq u(0,x)
\end{eqnarray*}
This inequality is true for all $h \in {\mathcal{H}}(T)$ so we have,
\begin{eqnarray*}
\sup_{h \in {\mathcal{H}}(T)} J(f,x,h,\gamma,T)  \leq u(0,x)
\end{eqnarray*}
Hence we have,
\begin{eqnarray}\label{5.14}
\inf_{\gamma \in \Gamma(T)}\sup_{h \in {\mathcal{H}}(T)} J(f,x,h,{\gamma},T)  \leq \sup_{h \in {\mathcal{H}}(T)} J(f,x,h,\gamma,T)  \leq u(0,x)
\end{eqnarray}
Likewise, setting $t=0$ and taking condition  (2) of Proposition 4.1 into account for $h=\hat{h}$, and for any $\gamma \in \Gamma(T)$ we see that
\begin{eqnarray*}
J(f,x,\hat{h},{\gamma},T) \geq u(0,x)
\end{eqnarray*}
This inequality is true for all $h \in {\mathcal{H}}(T)$ so:
\begin{eqnarray*}
\inf_{\gamma \in \Gamma(T)} J(f,x,\hat{h},{\gamma},T) \geq u(0,x)
\end{eqnarray*}
Hence we have,
\begin{eqnarray}\label{5.15}
\sup_{h \in {\mathcal{H}}(T)} \inf_{\gamma \in \Gamma(T)} J(f,x,{h},{\gamma},T)\geq \inf_{\gamma \in \Gamma(T)} J(f,x,\hat{h},{\gamma},T) \geq u(0,x)
\end{eqnarray}
Hence from (\ref{5.14}) and (\ref{5.15}) we have,
\begin{eqnarray}\label{5.16}
\sup_{h \in {\mathcal{H}}(T)} \inf_{\gamma \in \Gamma(T)} J(f,x,{h},{\gamma},T)\geq u(0,x)\geq \inf_{\gamma \in \Gamma(T)}\sup_{h \in {\mathcal{H}}(T)} J(f,x,h,{\gamma},T)
\end{eqnarray}
Moreover, setting $t=0$ and taking condition  (3) of Proposition 4.1 into account for $h=\hat{h},\gamma=\hat{\gamma}$ (such that $\hat{h} \in \mathcal{H}(T)$ and $\hat{\gamma} \in \Gamma(T)$) we see that
\begin{eqnarray}\label{5.17}
J(f,x,\hat{h},\hat{\gamma},T)= u(0,x)
\end{eqnarray}
It is always true that
\begin{eqnarray}\label{5.18}
\sup_{h \in \mathcal{H}(T)}(\inf_{\gamma \in  {{\Gamma}}(T)}J(f,x,h,{\gamma},T))\leq \inf_{\gamma \in  {{\Gamma}}(T)}(\sup_{h \in \mathcal{H}(T)} J(f,x,h,{\gamma},T))
\end{eqnarray}
Hence combining (\ref{5.16}) and (\ref{5.18}) we deduce the final conclusion  that the game ({\bf P1}) has a value and is $u(0,x)$. $\qed$\\

\section{Conclusion}
\indent In this article we provide a two player zero sum stochastic differential game in the context of the risk-sensitive benchmark asset management problem. We obtain an explicit expression for the optimal strategies for both the players. Future work could be directed towards considering a game theoretic benchmark problem with infinite horizon risk sensitive  criterion.

\end{document}